\documentclass[11pt]{article}
\usepackage{amsmath,amsfonts}
\usepackage{tensind}
\usepackage{graphicx}
\tensordelimiter{?}

\title{The Hyperbolic Heisenberg and Sigma Models in (1+1)-dimensions}
\author{A.E. Winn\footnote{E-mail: A.E.Winn@durham.ac.uk} \\
Department of Mathematical Sciences,\\
University of Durham, Durham DH1 3LE}
\date{}

\begin{document}
\maketitle

\begin{abstract}
Hyperbolic versions of the integrable (1+1)-dimensional Heisenberg
Ferromagnet and sigma models are discussed in the context of
topological solutions classifiable by an integer `winding number'. Some
explicit solutions are presented and the existence of certain classes of such
winding solutions examined.
\end{abstract}
\section{Introduction}
The classical Heisenberg model and the non-linear sigma model are well known
(1+1)-dimensional integrable systems which have been studied for both compact
and non-compact target spaces, c.f. [1-12].
The models are both mathematically and physically relevant for example,
(pseudo-)Riemannian manifolds are of great geometrical interest and a 
classical solution of such models is simply a harmonic map into such a 
space \cite{AntP2}.  Further, the Heisenberg model has been shown to be 
gauge equivalent to non-linear Schr\"{o}dinger equations of attractive
(in the compact case) \cite{FT} and repulsive (for a non-compact field space)
\cite{Kundu1} types, whilst a reduction of the $O(2,1)$ non-compact
sigma model is gauge equivalent to the Liouville equation \cite{DeVS}.
Physically, the Heisenberg model describes, in the compact case for example,
the classical spin $\vec{S}$ distributed along the line, i.e. the 
1-dimensional continuous ferromagnet. And the sigma models on non-compact
manifolds have arisen in relativistic string theory, see for example
\cite{DeVS}, and gravitation \cite{Mazur}.

In what follows, the discussion of solutions classified by an integer 
winding number of the hyperbolic (1+1)-dimensional Heisenberg (HHM) and sigma 
(HSM) models begun in \cite{US} is continued.  The field 
$\vec{\psi}(t,x) = (\psi^{1},\psi^{2},\psi^{3})$ takes values on the 
hyperboloid of one sheet $H^{2}$ in $\mathbb{R}^{2+1}$ and satisfies the
constraint
\begin{equation}
\eta_{ab} \, \psi^{a} \psi^{b} = 1, \label{eq:constraint1}
\end{equation}
with $\eta_{ab} = \mbox{diag }(1,1,-1)$. Also, $t \, \in \,
\mathbb{R}$ denotes time and $x \,\in \, X$ is the space variable such
that either $X = \mathbb{R}$ and the boundary condition $\vec{\psi}
(t,\infty) = \vec{\psi} (t,-\infty)$ is imposed, or $X = S^{1}$ so
that $\vec{\psi}$ is periodic in $x$.  With the metric on the
hyperboloid (\ref{eq:constraint1}) taken to be that induced by
$\eta_{ab}$, the manifold is a symmetric space $SO(2,1)/SO(1,1)$. Its
fundamental group then equals the group of integers $\mathbb{Z}$, so
that for each fixed $t$, $\vec{\psi}$ is a continuous mapping from a
circle into $H^{2}$ and therefore has a winding number  $N \, \in \,
\mathbb{Z}$ which is constant in $t$.  The field can be  visualized as
a closed string wrapped around the hyperboloid and evolving in time.

The Heisenberg model arises from the Hamiltonian density
\begin{equation}
\mathcal{H} = \frac{1}{2} \psi _{x}^{a} \psi _{x}^{b} \eta_{ab}
\label{eq:HMHamden}
\end{equation}
(where the subscript $x$ denotes differentiation with respect to $x$),
with Poisson brackets
\begin{equation}
\{ \psi ^{a}(x), \psi ^{b}(y) \} = -\delta (x - y) \epsilon ^{abc}
                                    \psi ^{d}(x) \eta _{cd}
                                    \label{eq:HMPoisson}
\end{equation}                                    
and the Hamiltonian energy for the system is given by
\begin{equation}
H = \frac{1}{2} \int_{X} \mathcal{H} \, dx, \label{HamE}
\end{equation}
which we expect to be non-positive-definite due to the indefinite metric on 
the field space. 
The equations of motion for the system are given by
\begin{equation}
\frac{\partial}{\partial t} \psi _{a}(x) = - ?{\epsilon}_{a}^{bc}? 
                                           \frac{\partial ^{2} \psi _{b}}
                                           {\partial x^{2}} \; \psi _{c}.
                                           \label{eq:HMmotion}
\end{equation}

The sigma model is defined by the Lagrangian density
\begin{equation}
\mathcal{L} = \frac{1}{2} \eta^{\mu\nu} \psi^{a}_{\mu} \psi^{b}_{\nu} 
              \eta_{ab}  \label{eq:Lagden}
\end{equation}
from which the equations of motion
\begin{equation}
(\psi^{a}_{\mu\nu} + \psi^{a}(\psi^{b}_{\mu} \psi^{c}_{\nu}) \eta_{bc})
       \eta^{\mu\nu} = 0   \label{eq:Sigmamotion}
\end{equation}       
arise, and where here and in what follows $\eta^{\mu\nu} =
\mbox{diag}(1,-1)$ and $x^{\mu} = (t,x)$.  Both models in
(1+1)-dimensions are integrable, in the sense that a suitable Lax pair
exists for each model which can be obtained by analytic continuation
from the analogous $S^{2}$ systems. 

\section{Travelling Waves for HHM}  \label{TWHM}

Parametrizing the hyperboloid in terms of `polar angles' $\theta
\mbox{ and } \phi$ so that 
\begin{equation}
\psi^{a} = (\cosh \theta \, \cos \phi, \cosh \theta \, \sin \phi, \sinh 
            \theta), \label{H2par1}
\end{equation}
equation (\ref{eq:HMmotion}) is then equivalent to 
\begin{eqnarray}
\theta _{t}  & = & 2(\sinh \theta) \theta _{x} \phi _{x} + (\cosh \theta )
                    \phi _{xx}  \label{eq:HMthetamotion}  \\
\phi _{t}    & = & (\mbox{sech} \theta ) \theta _{xx} + (\sinh \theta )  
                    \phi_{x}^{2}.  \label{eq:HMphimotion}
\end{eqnarray}
In \cite{US} some simple static and time dependent `winding' solutions
were shown to exist for this model, and here we extend the search for 
such winding solutions to those of travelling-wave type.  To facilitate this,
we introduce the functions $f(\xi) = \theta (t,x) \mbox{ and } g(\xi) + ct = 
\phi (t,x)$ with characteristic variable $\xi = x - vt$, where $v$ is the 
speed of the wave and interactions with an external magnetic field are 
introduced through the constant $c$.  Substitution into (\ref{eq:HMmotion})
results in 
\begin{equation}
\frac{dg}{d\xi} = (k - v \sinh \, f)\mbox{ sech }^{2}\, f  
      \label{eq:dgxi}
\end{equation}
where $k$ is constant and `prime' refers to differentiation with respect to 
$\xi$.  And this in turn leads to the equation for $f$:
\[ f'' = (v^{2} - k^{2}) \tanh f \mbox{ sech}^{2} \, f
                              + kv\mbox{ sech} \, f \, (1 - 2\mbox{ sech}^{2}
                              \, f) + c \, \cosh \, f. 
\]                              
This has first integral
\[ {f'}^{2} = \frac{(v^{2} - k^{2}) \sinh ^{2} \,f
    - 2kv \, \sinh \, f + 2c \, \sinh \, f \, \cosh ^{2} \, f + 2Q
    \cosh ^{2} \, f}{\cosh ^{2} \, f}
\]
with $Q$ constant which can be simplified by letting $p = \sinh 
\, f$. Substitution into the above then results in
\begin{equation}
\frac{{p'}^{2}}{2} = cp^{3} + \frac{p^{2}}{2} \, (v^{2} - k^{2} + 2Q) + 
                   p \, (c - kv) + Q  \, = P(p)  \label{eq:P(p)}
\end{equation}
with the 4 real parameters $k, v, c, Q$.
For a solution of winding type, we require first that $f(\xi)$ be in some 
sense periodic in $x$ (depending on the boundary conditions) and that the
integral over $X$ of (\ref{eq:dgxi}) be a finite integer multiple of $2\pi$
i.e. that, in terms of $p = \sinh f$
\[ \triangle \phi = \int_{X} \frac{dg}{d\xi} \, dx = \int_{X} 
   \frac{k - vp}{1 + p^{2}} \, dx = 2\pi N 
\]   
with $N \,\in \, \mathbb{Z}$. Taking as a first case, $c = 0$, i.e. in the
absence of an external field, one has 
\begin{equation}
\frac{{p'}^{2}}{2} = \frac{\alpha}{2} p^{2} - kvp + Q = P(p)
\label{eq:czeroP}
\end{equation}
where $\alpha = v^{2} - k^{2} + 2Q$.  
For a winding solution the restrictions $\alpha < 0 \mbox{ and }
Q > \frac{k^{2}v^{2}}{2\alpha}$ are necessary so that P(p) is positive and
bounded between two real zeros.  Reorganising the constants and integrating
results in the solution
\begin{equation}
p(\xi) = p_{0} + \frac{\sqrt{2A}}{N} \sin \bigl[ N(\xi - \xi)_{0} \bigr]
       \label{soln:sinNHM}
\end{equation}
with $A > 0, \, N \, \in \, \mathbb{Z}, p_{0} \mbox{ and } \xi_{0}$ constant,
so that $p(\xi) = \sinh\, f$ is $2\pi$ periodic. (We note that this is a
travelling-wave version of the static solution (13) in \cite{US}).  That this
is indeed a winding solution can be seen as follows: taking for simplicity
$N = 1$ (since it is a fairly simple generalization to $\mid N \mid > 1$),
scaling $\alpha = -1$ thereby eliminating $Q$, and with
$\xi_{0} = 0$, in terms of our original constants $k, v$, the solution
(\ref{soln:sinNHM}) can be rewritten
\begin{equation}
p(\xi) = -kv + \sqrt{(v^{2} + 1)(k^{2} - 1)}\, \sin (\xi)
\label{eq:c0p1}
\end{equation}
with period $2\pi$.  We now have two parameters: 
$v \, \in \, \mathbb{R} \,\mbox{ and } k \, \geq \, 1$ so that 
\begin{description}
\item[(i)] if $v = 0$, then $\triangle \phi = 2\pi \; \forall \, k
\, \geq \, 1$,  c.f. analysis in  {\textbf{(iii)}} below,  
\item[(ii)] if $k = 1$, then $\triangle \phi = 2\pi  \; \forall \,
v$ (where this winding solution has $p \equiv -v)$,  
\item[(iii)] if $v \neq 0 \, \mbox{ and }k \, > \, 1$ (\ref{eq:dgxi}) 
reduces to 
\[ \frac{dg}{d\xi} = \frac{1}{2} \left[ \frac{(v + ik)}{(kv + i) - \gamma
                     \sin \xi} + \frac{(v - ik)}{(kv - i) - \gamma \sin \xi}
                     \right].
\]                     
This integrates to
\begin{equation}
g(\xi) =  \mbox{Re } \Biggl[ -\frac{i}{2} \left[ \ln \left\{ \frac{\tan 
          \frac{\xi}{2} -
          \frac{(i\gamma + ik - v)}{(1 + ikv)}}{\tan \frac{\xi}{2} -
          \frac{(i\gamma - ik + v)}{(1 + ikv)}} \right\}
          - \ln \left\{ \frac{ \tan \frac{\xi}{2} + 
          \frac{i\gamma + ik + v}{(1 - ikv)}}{\tan \frac{\xi}{2} + 
          \frac{i\gamma - ik - v}{ (1 - ikv)}} \right\} \right] + \Lambda
          \Biggr]
\end{equation}
where $\Lambda$ is a constant of integration and may be simplified to
\begin{equation}
g(\xi) = \mbox{Re } \left[ -i \ln \left\{ \frac{\tan \frac{\xi}{2} +
         \frac{a}{b}}{\tan \frac{\xi}{2} + \frac{c}{b}} \right\} + \Lambda  
         \right], \label{soln:Re1}
\end{equation}  
where $a = -i \gamma -ik + v, \, b = 1 + ikv, \, c = -i \gamma + ik - v$.
As it stands, it is not easy to see that this solution is of winding
type; however, it can be shown that (\ref{soln:Re1}) reduces to the form 
\begin{equation}
g(\xi) = \mbox{Re } \left[ -i \ln \left\{ \frac{ \tan \frac{\xi}{2} +
         \Omega}{-\tan \frac{\xi}{2} - \Omega ^{-1}} \right\} \right]
         \label{soln:Re2}
\end{equation}
where $\Omega = \frac{v - ik - i\gamma}{1 + ikv} \left( = \frac{a}{b}
\right)$. And this {\em{is}} of winding type, as can be seen by analysing the
behaviour of the function
\[ \Xi (x) = \frac{ \tan \frac{\xi}{2} + \Omega}{- \tan \frac{\xi}{2}
              - \Omega ^{-1}}
\]
as $x$ goes from $- \pi \mbox{ to } \pi$. Briefly, this is as follows:
$\Xi (x)$ is a continuous function of $x$
which is never zero (for $\Xi$ to be zero requires $x \, \in \,
\mathbb{C}$)
and as $x \longrightarrow \pm \pi ,  \, \Xi \longrightarrow -1$ so that
$\Xi$ takes the form of a loop, beginning and ending at $-1$.  Now taking
a branch cut from the origin of $\Xi$ space along the negative real axis,
it remains to show that $\Xi$ wraps once around the origin so that   
$\triangle \phi  = \pm 2\pi$.
For this to be the case, $\Xi$ must pass through the positive real axis
once and only once so that $\Xi = \overline{\Xi}$ must have a {\em{unique}},
real solution $x$. This is indeed the case since $\Xi = \overline{\Xi}$ 
requires that
\begin{eqnarray*}
\frac{ \tan \frac{\xi}{2} + \Omega}{ - \tan \frac{\xi}{2} - \Omega ^{-1}}
& = & \frac{ \tan \frac{\xi}{2} + \overline{\Omega}}{- \tan \frac{\xi}{2}
      - \overline{\Omega} ^{-1}} \\
\Longrightarrow \hspace{.2in} \tan \frac{\xi}{2}
& = & \frac{ \overline{\Omega} \Omega ^{-1} - \Omega
      \overline{\Omega}^{-1}}{\Omega + \overline{\Omega}^{-1} - \Omega
^{-1}
      - \overline{\Omega}}.
\end{eqnarray*}
so that
\[ \xi = 2 \tan ^{-1} \left[ \frac{ -2 (v - k^{2}v - \gamma kv)
        (1 + k^{2}v^{2})}{(v - k^{2}v - \gamma kv )^{2} + (k + \gamma +
        kv^{2} )^{2} + (1 + k^{2}v^{2})^{2}} \right]   
\]
giving a unique and real solution $x$ for
each $k, v \, \in \mathbb{R}$. Then
$\phi = \pm 2 \pi$, i.e. the solution (\ref{soln:Re2}) is a travelling
wave winding solution for the HHM where the space $X = S^{1}$ and we 
have taken $c = 0$. Hence  (\ref{soln:Re1}) is also such a solution.
\end{description}

In our next example, we take the constant $c = 1$, in an
attempt to find a more general solution. With the
substitutions $\alpha = \frac{1}{6}(k^{2} - v^{2} + 2Q), \; \triangle
= \alpha ^{2} - \frac{1}{3} (1 - kv)$, equation (\ref{eq:P(p)}) is
transformed to
\begin{equation}
P(p) = \frac{{p'}^{2}}{2} =  
      p^{3} - 3p^{2} \alpha + 3(\alpha^{2} - \triangle)p + Q
      \label{eq:c1P(p)}
\end{equation}
retaining the 3 real parameters $k, v, Q$. And 
for a solution to be a travelling wave of winding type, $P(p)$ should 
have  three real
zeros, so that there exists a positive bounded region in which the
solution can oscillate.  This imposes the restriction that 
$\triangle > 0$. We shall presently examine how the limits of
this behaviour are manifested for the spaces $X = S^{1} \mbox{ and } X
= \mathbb{R}$; however, to facilitate this we first note that if $P(p)$,
(\ref{eq:c1P(p)}), has three real and distinct zeros $p_{1} > p_{2} > p_{3}$
then the following solution can be found in terms of elliptic funtions:
\begin{equation}
p(\xi) = p_{2} - (p_{2} - p_{3}) \mbox{cn}^{2} \left[ (\xi - \xi_{3})
\left\{ \frac{p_{1} - p_{3}}{2} \right\}^{\frac{1}{2}} \bigl| m \bigr. 
\right] \label{soln:DJ}
\end{equation}
where the parameter $m = \frac{p_{2} - p_{3}}{p_{1} - p_{3}}$.
Note that this is also a travelling wave solution of the KdV
equation \cite{DJ}, and given values of $p_{1}, p_{2}, p_{3}$, the shape
of the cnoidal wave can be obtained either from tables of Jacobian
elliptic functions or by direct computation.

To determine whether or not
the solution is actually of winding type for our model if $p(\xi)$ is given 
by (\ref{soln:DJ}),  
it is no simple matter (and may not be possible at all) to integrate 
(\ref{eq:dgxi}). So rather than attempt this we will investigate what
happens  in the limiting cases, i.e. where the parameter 
$m = 0 \mbox{ or } 1$. When $m \equiv  0$ one has $p_{2} = p_{3}$ 
and $P(p)$ has one repeated and one single zero with no positive
bounded region. 
The solution is then $p = p_{3} = \sinh f$ and in this case, $\frac{dg}{d\xi} 
= \frac{k - vp_{3}}{(1 + p_{3}^{2})} =  N \mbox{ constant }$; so that if 
$N \,\in \, \mathbb{Z}$ we recover the winding solution
\[ \phi = Nx + t
\]
of \cite{US}.
A further limiting case occurs when $P(p)$ has one repeated and one
single zero where the region between them is positive and bounded.  If
$p_{1}$ is the double zero and $p_{3}$ the single, then
this picture corresponds to
a travelling wave on $X = \mathbb{R}$ where the solution has a minimum
at $p= p_{3}$ (since $P'(p_{3}) > 0 $) and reaches $p = p_{1}$ as $\xi
\longrightarrow \pm \infty$. Now referring back to solution
(\ref{soln:DJ}) and considering the limit where the parameter $m =
\frac{p_{2} - p_{3}}{p_{1} - p_{3}} = 1$, one can see
that $p_{2}$ must equal $p_{1}$ which is
exactly the situation described.  And since 
$\mbox{ cn }(u|m) = \mbox{ sech }(u) \mbox{ for } m = 1$, 
(\ref{soln:DJ}) reduces to 
\begin{equation}
p(\xi) = p_{1} - (p_{1} - p_{3}) \mbox{ sech }^{2} \Bigl[ (\xi - \xi
_{0}) \Bigl\{ \frac{p_{1} - p_{3}}{2} \Bigr\} ^{\frac{1}{2}} \Bigr].
\label{soln:HMm1} 
\end{equation}
It is possible to find explicit expressions for 
$p_{1}, \, p_{3}$ satisfying our equation (\ref{eq:c1P(p)}) and giving the
correct picture, however, it turns out that this is not a winding solution 
for our Heisenberg model; on substitution of (\ref{soln:HMm1}) into 
(\ref{eq:dgxi}) one finds, for example in the one-wind case, that
\begin{equation}\begin{split}
\frac{dg}{d\xi}&= \frac{k - v\gamma}{(\gamma^{2} - 1)} + 
                  \Omega \Bigl[ \frac{1}{\lambda - \cosh X} +
                  \frac{1}{\lambda + \cosh X} \Bigr]  \label{HMm1dg} \\
               & + \overline{\Omega} \Bigl[ \frac{1}{\overline{\lambda}
                  - \cosh X} + \frac{1}{\overline{\lambda} + \cosh X}
                  \Bigr]
\end{split}\end{equation}
where $X = (\xi - \xi_{0}), \, \gamma = \frac{1}{6}(k^{2} - v^{2} + 2Q + 4),
\, \Omega = \frac{(v - ik)\sqrt{\gamma + i}}{2\sqrt{2} ( \gamma +
i )^{2}} \mbox{ and } \lambda = \bigl[ \frac{2}{\gamma + i }\bigr]
^{\frac{1}{2}} $.  Integrated over $X = \mathbb{R}$, (\ref{HMm1dg})
diverges so that $\phi$ wraps infinitely many times around the
hyperboloid $H^{2}$. And this solution  is therefore not of
winding type since we require that the winding number $N \, \in \,
\mathbb{Z}$ be finite.

We note further that this behaviour also manifests itself in
the corresponding Hamiltonian density. For example,
take the simple case $p_{1} = 0$, so that $Q = 0$ and for the one-wind
solution, $\sqrt{\triangle} = -\frac{2}{3}, \, v^{2} = 2 + \sqrt{5}
\mbox{ and } p_{3} = -2$.
The Hamiltonian density
is then given by 
\[ \mathcal{H} = \frac{\bigl[ 1 + 2(2 + \sqrt{5}) \mbox{sech}^{2}(X)
                 \bigr]^{2} - 16(\sqrt{2 +
                 \sqrt{5}})\mbox{sech}^{4}(X)\tanh ^{2}(X)}{(\sqrt{2 +
                 \sqrt{5}}) (1 + 4\mbox{ sech}^{4}(X))},
\]
the profile of which
\begin{figure}
\begin{center}
\includegraphics[angle=-90,width=4cm]{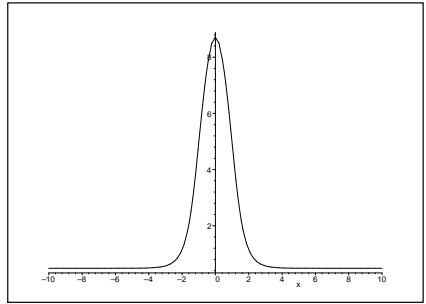}
\caption{Hamiltonian Density for the one-wind solution when $c = 1$ in
                 the limit $m= 1$ with $p_{0} = 0$}
\label{fig:HE}
\end{center}
\end{figure}
is shown in figure (\ref{fig:HE}).  One can
see that as  $X = \xi - \xi_{0} \longrightarrow \pm \infty,\,
\mathcal{H}$ does not tend to zero hence, the Hamiltonian energy is infinite. 
Since the profile of the density
$\mathcal{H}$ is, however, that of an attractive uniform lump, some
renormalization procedure might be utilized to remove this infinity
but we will not pursue this possibility here.

In \cite{US} we stated that in fact, no travelling wave solutions of 
winding type exist with $X = \mathbb{R}$ for this model and we now show this 
to be the case. For such solutions to exist requires
\begin{description}
\item[(i)] the solution $p(\xi)$ be real and bounded so $P(p) = {p'}^{2} 
           \geq 0$ (from the form of $P(p)$ (\ref{eq:P(p)}) one can see 
           that $p(\xi)$ will vary monotonically until $p'(\xi)$ vanishes);
\item[(ii)] $P(p)$ have two real zeros; one repeated and one simple ($p_{0}$
           and $p_{1}$ respectively), where the solution $p(\xi)$ has a 
           minimum at $p = p_{1}$ and $p(x) \longrightarrow p_{0} 
           \mbox{ as } x \longrightarrow \pm \infty$.
\end{description}
Recalling further that 
\[ \triangle \phi = \int_{X} \frac{dg}{d\xi} \, dx = \int_{-\infty}^{\infty}
   \frac{k - vp(\xi)}{1 + p^{2}(\xi)} \, dx
\]
we must further demand that
\begin{description}
\item[(iii)] $p(\xi) \longrightarrow \frac{k}{v} \mbox{ as } x \longrightarrow
           \pm \infty$ where $\frac{k}{v}$ is the zero of $\frac{dg}{d\xi}$,
           so that $p_{0} = \frac{k}{v}$ ({\textbf{(ii)}} above), and 
\item[(iv)] $p'(\xi) \longrightarrow M$ const. as $x \longrightarrow 
            \pm \infty$ and in fact on $X = \mathbb{R}$,
            {\textbf{(iii)}} and {\textbf{(i)}} above
            imply $M = 0$, so that $P(p = \frac{k}{v}) = 0$. 
\end{description}
We note first that with $p_{0} = \frac{k}{v}, \; P(p)$ must have the
form
\begin{equation}
P(p) = \left( p - \frac{k}{v} \right) ^{2} \left( \alpha p - \beta
\right)  \label{comparepzero}
\end{equation}
for some $\alpha, \, \beta \, \in \, \mathbb{R}$.  To show that the
conditions {\textbf{(i)}}-{\textbf{(iv)}} cannot be satisfied for a
solution $p = \sinh f$ of the HHM we need only now use condition
{\textbf{(iv)}}, i.e. $P\left( p = \frac{k}{v} \right) = 0$ in
equation (\ref{eq:P(p)}). This substitution results in the equation
\begin{equation}
(k^{2} - v^{2})(2Q + 2ck - k^{2}v^{2}) = 0  \label{kovervzero}
\end{equation}
presenting the following possibilities:
\begin{equation}
\begin{array}{ll}
(\mbox{a}) & \; v \; = \; \pm ik \\
(\mbox{b}) & \; v \; = \; c \; = \; 0 \\
(\mbox{c}) & \; v \; = \; k \; = \; 0 \\
(\mbox{d}) & \; k \; = \; Q \; = \; 0 \\
(\mbox{e}) & \; Q \; = \; \frac{k}{2v}(kv - 2c) \\
(\mbox{f}) & \; c \; = \; \frac{v}{2k}(k^{2} - 2Q). 
\end{array} \nonumber
\end{equation}
The first three can be immediately discounted since we require real,
non-zero $v$ and case (d) results in $P(p) = \frac{v^{2}}{2} p^{2}$,
thereby contradicting {\textbf{(ii)}} above.  Substituting (e) into
(\ref{eq:P(p)}) and comparing the resultant equation with
(\ref{comparepzero}) allows only $v = \pm ik$ so that we are left with
case (f). Following the same procedure  as with (e) one finds that
either $v = \pm ik$ or $k^{2} = 2Q$.  In the latter instance, this
imposes that $P(p) = \frac{1}{2}\left( p - \frac{k}{v} \right) ^{2}$,
again contradicting {\textbf{(ii)}}.
All possibilities have then been exhausted for values of the
parameters satisfying (\ref{comparepzero}) so that no $p = \sinh f$
can exist satisfying all of the conditions
{\textbf{(i)}}-{\textbf{(iv)}}. Hence, whilst they do exist for the
space $X = S^{1}$, we have shown that there are no travelling wave
solutions of winding type for this model when $X$ is the real line
$\mathbb{R}$. 

\section{The Hyperbolic Sigma Model (HSM)}

With the parametrization (\ref{H2par1}) of the hyperboloid, the equations
of motion for the  HSM are 
\begin{eqnarray}
\theta_{tt} - \theta_{xx}         & = & - \cosh \theta \, \sinh \theta \, 
                                        (\phi_{t}^{2} - \phi_{x}^{2} ) 
                                        \label{eq:Sigmathetamotion} \\
( \phi_{t} \cosh^{2} \theta )_{t} & = & ( \phi_{x} \cosh^{2} \theta )_{x}.
                                        \label{eq:Sigmaphimotion}
\end{eqnarray}
The  energy is given by $E = \int_{X} \epsilon \, dx$ where
the energy density $\epsilon$ is given here by
\begin{equation} 
\epsilon = \cosh ^{2} \theta \, (\phi_{t}^{2} + \phi_{x}^{2}) -
              (\theta_{t}^{2} + \theta_{x}^{2}). \label{edensityHSM}
\end{equation}              
Note that this is indefinite due to the signature of the 
metric on $H^{2}$.

If one tries to find a solution with $\phi = Nx$ it turns out that
$\theta$ is a function of $t$ only, in which case the equations admit
the solution 
\begin{equation}
\tanh \theta = \mbox{ sn }\bigl[ \rho (t - t_{0}) | m \bigr], 
                \label{soln:tanhtheta}
\end{equation}
with $m = 1 - \frac{N^{2}}{\rho^{2}}, \, (\rho, t_{0} \mbox{ constants }),
\, \rho >  | N |$.  This is almost identical to a solution of the equations
in the positive definite version \cite{US}.
Observing that $\tanh ^{-1} u \longrightarrow \infty \mbox{ as } u 
\longrightarrow 1$, one can see that $\theta(t) \longrightarrow \infty 
\mbox{ as  sn}\,(\rho(t - t_{0}) | m ) \longrightarrow 1$
i.e. $\theta(t)$ reaches infinity in a finite time.  
The energy for this solution is given by 
$E = 2\pi \bigl[ A(\theta) - B(\theta) \bigr] \,$ where $A(\theta) = 
N^{2} \cosh ^{2} \theta$ and $B(\theta) = N^{2} \sinh ^{2} \theta + \rho^{2}$.
Each of these positive functions reaches infinity in 
a finite time (this instability being explicit in the solutions), however, it
is interesting to note that the total energy is a conserved 
(negative) quantity. The same calculation shows that in the positive definite 
case of \cite{US}, this quantity is 
greater than zero and indeed greater than the lower bound for the energy.

In our search for solutions of the system, as with the previous section, we
turn our attention to travelling waves for the HSM.  The equations
of motion in this case, with $\theta = f(\xi), \, \phi = g(\xi) + ct$
are then
\begin{eqnarray}
f'' (v^{2} - 1) & = & - \left( {g'}^{2} (v^{2} - 1) - 2vcg' + c^{2} \right)
                       \cosh f \, \sinh f   \label{eq:Sf''} \\
g'' (v^{2} - 1) & = & -2 \left( f'g' (v^{2} - 1) - cvf' \right) \tanh f.  
                      \label{eq:Sg''}
\end{eqnarray}
As a first case let us again take $c = 0$ so that
\begin{eqnarray}
f'' & = & - {g'}^{2} \cosh f \, \sinh f   \label{HSMf2} \\
g'' & = & -2 f' g' \tanh f. \label{HSMg2}
\end{eqnarray}
Noting then that
\begin{eqnarray*}
\frac{d}{d\xi} (g' \cosh ^{2} f)  & = & 2 f'g' \cosh f \, \sinh f
                                           + g''\cosh ^{2} f  \\
                                     & = & 0
\end{eqnarray*}
by (\ref{HSMg2}), 
we have the following first integrals for $g \mbox{ and } f$:
\begin{eqnarray}
g' & = & B \mbox{ sech}^{2} f  \label{eq:Sg'}  \\
f' & = & \sqrt{(B^{2} - N^{2})  - B^{2} \tanh ^{2} f} \label{eq:Sf'}
\end{eqnarray}
where $B \mbox{ and } N$ are constants.  
An analogous solution to (\ref{soln:sinNHM}) derives from these equations,
namely
\begin{equation}
p(\xi) = \sqrt{ \frac{B^{2}}{N^{2}} - 1} \sin \bigl[ N(\xi - 
                   \xi_{0}) \bigr]
\end{equation}
where again $p = \sinh f,  \mbox{ and } \xi_{0}$ is a constant.  This
is a winding solution  for HSM as can be shown in much the same way as
we did for the solution (\ref{soln:sinNHM}) of HHM. 

If $c \neq 0$, the following equations arise: 
\begin{eqnarray}
P(p) = {p'}^{2} & = & \frac{c^{2}}{(v^{2} - 1)^{2}} \Bigl[ 
               p^{4} + p^{2} 
               \frac{(6c^{2} + 4Q(v^{2} - 1)^{2})}{4c^{2}} \nonumber \\  
               &  & + \frac{2c^{2} + 4Q(v^{2} - 1)^{2} + \bigl
               ( 4R(v^{2} - 1)^{2} - cv \bigr)^{2}}{4c^{2}}
               \label{HSMpprime} \\
g' & = & \frac{cv + 2R(v^{2} - 1) + 2cvp^{2}}{(1 + p^{2})}   
\end{eqnarray}
where $R \mbox{ and } Q$ are constants.  And for $\sinh f$ to be periodic
$P(p)$ must take the form
\begin{equation} 
P(p) = (p^{2} - J^{2})(p^{2} - K^{2}) \label{P(p)JK}
\end{equation}
for $J, K \, \in \mathbb{R}$
so that $P(p)$ is positive and bounded between the two zeros $p = \pm J$.
Comparing coefficients results in the following two equations for $J 
\mbox{ and } K$:
\begin{eqnarray}
-(J^{2} + K^{2}) & = & \frac{6c^{2} + 4Q(v^{2} - 1)^{2}}{4c^{2}}  \nonumber \\
J^{2}K^{2} & = & \frac{2c^{2} + 4Q(v^{2} - 1)^{2} + \bigl(4R(v^{2} - 1)^{2}
                 - cv \bigr)^{2}}{4c^{2}}.  \label{eq:JK}
\end{eqnarray}
which, putting $q = Q(v^{2} - 1)^{2}, \, r = R(v^{2} - 1)^{2}, 4\rho
=4r - cv$ and without
loss of generality $c = 1$, imposes the conditions $q < -\frac{3}{2}, \;
2q > -(8\rho^{2} + 1)$ and $1 - 2q \geq 8 |\rho |$.  Such $q,
\rho$ can be found so that the solution for $p(\xi)$ in terms of 
elliptic functions is then given after a shift in $x$, by 
\begin{equation}
p(\xi) = J \mbox{ sn }\left[ \frac{cK(\xi - \xi_{0})}{(v^{2} - 1)} | m 
         \right]  \label{ellHSM}
\end{equation}
with $m = \frac{J^{2}}{K^{2}}$.  Examining the limiting cases of this solution
we find that with $m \equiv 0, \, p(\xi) \equiv 0$; and in the limit
$m = 1$, where we put $K = J = p_{0}$, 
\begin{equation}
p(\xi) = p_{0} \tanh \left[ \frac{cp_{0}(\xi - \xi_{0})}{(v^{2} - 1)} \right].
\end{equation}
In this case the polynomial $P(p)$ is positive and bounded between the two 
double zeros $p = \pm p_{0}$ so that the solution starts from $-p_{0}$ at 
$x = -\infty$ and travels to $p_{0}$ at $x = +\infty$, completing only a 
half period since it has effectively run out of $x$ space.  Since this is 
the only possible scenario for $X = \mathbb{R}$, this shows that there are no
travelling wave solutions of winding type for the model with $X = \mathbb{R}$.

If $P(p)$ takes the form (\ref{P(p)JK}), for winding solutions
we would, of course, expect $\triangle \phi = 2\pi N$ as usual and if this 
were the case, in the limiting case just described where $m = 1$ one would
expect $\triangle \phi$ to amount to half that since one
half period is traversed. However, this is not the case since
\[
\frac{dg}{d\xi} = \frac{cv}{(v^{2} - 1)} + \frac{2R(v^{2} - 1) - 
                  cv}{2(v^{2} - 1)\bigl[ 1 + p_{0}^{2} \tanh ^{2} X \bigr]}
\]
where $X = \pm \frac{cp_{0}(\xi - \xi_{0})}{(v^{2} - 1)}$ and the integral
of this over $\mathbb{R}$ clearly diverges.  So 
whilst this does not constitute a proof, it at least gives us an indication
that the general solution (\ref{ellHSM}) may not be of winding type.

\section{Concluding remarks}
Integrable systems admitting time dependent topological solutions are
comparatively rare, the sine-Gordon being virtually the only well known
example. The models described here are such integrable models and we have shown
that whilst no travelling wave solutions classifiable by an integer
winding number exist for the hyperbolic Heisenberg model  where $X =
\mathbb{R}$ and that they are unlikely to exist for the sigma
model with this $X$, they do exist when $X = S^{1}$. In
fact, we have shown that there exist solutions with $X = S^{1}$ which,
whilst not being identical in the different models,  are of a common
form to both. It may therefore be interesting to attempt some form of
interpolation between the solutions and possibly the models themselves. 

\vspace{.3in}

\noindent{\textbf{Acknowledgements}}

\noindent This work was supported by an EPSRC Research
Studentship. AEW wishes to thank Prof. R.S. Ward for his
patience and guidance.

\end{document}